\documentclass[reqno]{amsart}

\usepackage{algorithm}
\usepackage{algorithmicx}
\usepackage{algpseudocode}
\usepackage{amsmath,amssymb,amsfonts}
\usepackage{amsthm}
\usepackage[title]{appendix}
\usepackage{booktabs}
\usepackage[autostyle]{csquotes}
\usepackage{caption}
\usepackage{hyperref}
\usepackage{graphicx}
\usepackage{listings}
\usepackage{principia}
\usepackage{manyfoot}
\usepackage{mathrsfs}
\usepackage{multirow}
\usepackage{pdfpages}
\usepackage{textcomp}
\usepackage{xcolor}

\newcommand{\keyw}[1]{\textit{Keywords---} #1}

\newtheorem{theorem}{Theorem}
\newtheorem{interpretation}{Interpretation}
\newtheorem*{interpretation*}{Interpretation}
\newtheorem{lemma}[theorem]{Lemma}
\newtheorem{proposition}[theorem]{Proposition}
\newtheorem{remark}{Remark}

\begin{document}

\title{On Chwistek's criticism of \emph{Principia}}
\author{Stephen Boyce}

\begin{abstract}
  This paper examines Chwistek's claim that with \emph{Principia}'s definition of a
  class ``Richard’s paradox can be formulated''. It is shown that the demonstration fails
  since it requires an incorrect elimination of a defined term and use of a faulty substitution rule, neither of which form part of the system of \emph{Principia}.
\end{abstract}

\maketitle
\begin{msc}
03-03
\end{msc} \\
\keyw{Principia Mathematica, paradox}

\section{Introduction}

Leon Chwistek gains a mention in the Introduction to the Second Edition of \emph{Principia Mathematica} in relation to the idea of dispensing with the axiom of reducibility
(\cite{pm1910v1}: xiv). My interest here however is with Chwistek's claim that with \emph{Principia}'s definition of a
class ``Richard’s paradox can be formulated'' (\cite{chwistek1922}: 338). The claim, if accepted, puts into disrepute any part of the First Edition involving set theory.
Due to an overwhelming lack of interest in \emph{Principia}, it is a little difficult to gauge contemporary logicians' perceptions of Chwistek's claim.
To pick a couple of commentators at random,
Linsky's \cite{linsky2011} assessment seems favourable as too does G\"odel's (\cite{godel1944}: 129).
Whilst G\"odel obviously does not qualify as a contemporary logician, his pronouncements on \emph{Principia} arguably carry a lot of weight with
contemporary scholars.

This then brings us to the question I wish to address. Does Chwistek's demonstration hold up?
\section{A reassessment of Chwistek's demonstration}
One accessible English-language version of the demonstration is as follows.
\begin{quote}
  Consider the theorem ($\pmast 20\pmcdot1$)
  \begin{equation*}
    f\{\pmcls{z}{\pmpf{\varphi}{z}}\}
    \pmdot \pmiff \pmdott
    \pmsome{\psi} \pmdott \pmpred{\psi}{x} \pmdot \pmiff_x \pmdot \pmpf{\varphi}{z} \pmandd f\{\pmpred{\psi}{\pmhat{z}}\}
  \end{equation*}
  and substitute the following function for $f\{\pmpred{\psi}{\pmhat{z}}\}$:
  \begin{equation*}
    \pmnot \{\pmcls{z}{ \pmpred{\chi}{z} }  = \pmpred{\psi}{\pmhat{z}} \}.
  \end{equation*}
  We get
  \begin{equation*} \label{equation_chwistek_1}
    \pmnot \{\pmcls{z}{ \pmpred{\chi}{z} }  = \pmcls{z}{ \pmpf{\phi}{z} }   \}
    \pmdot \pmiff \pmdott
    \pmsome{\psi} \pmdott \pmpred{\psi}{z}
    \pmdot \pmiff_z \pmdot
    \varphi z \pmand \pmnot \{\pmcls{z}{ \pmpred{\chi}{z} }  = \pmpred{\psi}{\pmhat{z}} \} \tag{1}
  \end{equation*}
  and by eliminating $\pmcls{z}{\varphi z}$. On the basis of a definition of a statement about class and the
  corresponding convention, we conclude
  \begin{align*}
    \pmnot \pmsome{\psi} \pmdott \pmpred{\psi}{z} \pmdot
    \pmiff_{\chi} \pmdot \pmpf{\varphi}{z} \pmdot
    \{\pmcls{z}{ \pmpred{\chi}{z} }  = \pmpred{\psi}{\pmhat{z}} \}
    \pmdott \pmiff \pmdott 
    \pmsome{\psi} \pmdott \pmpred{\psi}{z} \pmdot \pmiff_{\chi} \pmdot \varphi z \pmand &  \\ 
    \pmnot \{\pmcls{z}{ \pmpred{\chi}{z} }  = \pmpred{\psi}{\pmhat{z}} \}, \tag{2}  \label{equation_chwistek_2} &&
  \end{align*}
  which is evidently false. (\cite{chwistek1922}: 339 modified through addition of reference numbers for equations SB.)
\end{quote}
A careful assessment of this demonstration, supplemented by an examination of relevant related material \cite{chwistek1924}, shows that Chwistek's claim does not hold up.
\begin{proposition}\label{proposition_chwistek_0}
  Chwistek's demonstration fails to establish that Richard’s paradox can be formulated in \emph{Principia}.
  Chwistek's derivation of Equation \ref{equation_chwistek_2} is linked with an error in eliminating a defined term. 
  Chwistek's perception of the falsity of Equation \ref{equation_chwistek_2} is furthermore partly due to the use of a faulty substitution rule which would be prohibited by any
  reasonable statement of substitution rules for \emph{Principia}. 
\end{proposition}
\begin{proof}
  To obtain Equation \ref{equation_chwistek_1} from $\pmast 20\pmcdot1$, Chwistek uses the function $f$ such that:\footnote{Bearing in mind that
  ``$f\{\pmcls{z}{\pmpf{\varphi}{z}}\}$ is not a value of $fx$'' (\cite{pm1910v1}: 202) for reasons set out therein.}
  \begin{flalign*}
    f\{\pmcls{z}{\pmpf{\varphi}{z}}\} = \pmnot \{\pmcls{z}{ \pmpred{\chi}{z} }  = \pmcls{z}{\pmpf{\phi}{z}} \}\ \label{item_chwistek_item_1} \text{and} & \tag{3} \\
    f\{\pmpred{\psi}{\pmhat{z}}\} = \pmnot \{\pmcls{z}{ \pmpred{\chi}{z} }  = \pmpred{\psi}{\pmhat{z}} \}\ \label{item_chwistek_item_2} \text{etc.} && \tag{4}
  \end{flalign*}
  Chwistek's inference from Equation \ref{equation_chwistek_1} to Equation
  \ref{equation_chwistek_2} is justified as
  ``On the basis of a definition of a statement about class and the corresponding convention'' (\cite{chwistek1922}: 339). The ``definition of a statement about class''
  mentioned is a reference to $\pmast 20\pmcdot01$
  while the ``corresponding convention'' referred to is this:
  \begin{quote}
    The scope of the symbol ``$\pmcls{z}{\pmpf{\psi}{z}}$'' is the smallest proposition enclosed in dots or brackets in which  ``$\pmcls{z}{\pmpf{\psi}{z}}$'' occurs.
    [footnote:\textquote*{Cf Principia Vol.I p. 197 and 181.}] (\cite{chwistek1924} in \cite{Linsky2004}: 68)
  \end{quote}
  Chwistek's statement of a rule for determining the scope of the symbol ``$\pmcls{z}{\pmpf{\psi}{z}}$'' will generally give a correct result,
  that is, a determination of the scope in agreement with the broader practice of \emph{Principia}. The statement, as it stands, however obscures a subtle point that
  should be considered to properly evaluate Chwistek's demonstration.
  \begin{remark}\label{remark_context_0}
    Correct determination of the scope of the symbol ``$\pmcls{z}{\pmpf{\psi}{z}}$'' requires consideration of the relevant definitions, especially with respect to the
    context of use:
    \begin{enumerate}
    \item Where a class symbol occurs as part of a propositional function that is an extensional function of it's arguments, the above rule is guaranteed to give a correct result
      in virtue of the axiom of reducibility (\cite{pm1910v1}: 83);
    \item There are however contexts in which ``two functions may well be formally equivalent without being identical'' (\cite{pm1910v1}: 83). Russell gives the example of the two functions ``$x =$ Scott'' and
      ``$x =$ the author of Waverley''. If such contexts are in scope for a discussion then we cannot simply apply the above rule.
    \end{enumerate}
  \end{remark}
  To see how these issues affect Chwistek's argument, let's firstly retrace Chwistek's demonstration. We may start using Equation \ref{item_chwistek_item_1} and $\pmast 20\pmcdot1$
  obtaining:\footnote{I have adjusted some of \emph{Principia}'s notation to simplify comparison with Chwistek's presentation, e.g.
  ``$f\{\pmcls{z}{\pmpf{\psi}{z}}\}$'' is given as ``$f\{\pmcls{z}{\pmpf{\varphi}{z}}\}$'' etc.}
  \begin{flalign*} 
    \quad \pmSub{\pmast 20\pmcdot1}
          {\pmnot \{\pmcls{z}{ \pmpred{\chi}{z} }  = \pmcls{z}{\pmpf{\phi}{z}} \}}
          {f\{\pmcls{z}{\pmpf{\varphi}{z}}\}}\
          \pmthm \pmdottt
          \pmnot \{\pmcls{z}{ \pmpred{\chi}{z} }  = \pmcls{z}{ \pmpf{\phi}{z} }   \}
          \pmdot \pmiff \pmdott  & \\
          \pmsome{\psi} \pmdott \pmpred{\psi}{z} \pmdot
          \pmiff_{z} \pmdot  \pmpf{\phi}{z} 
          \pmdott
          \pmnot \{\pmcls{z}{ \pmpred{\chi}{z} }  = \pmpred{\psi}{\pmhat{z}} \} \tag{5} \label{equation_chwistek_4} &&
  \end{flalign*}
  Then, using a different function $f\{\pmcls{z}{\pmpf{\varphi}{z}}\} = \{\pmcls{z}{ \pmpred{\chi}{z} }  = \pmcls{z}{\pmpf{\phi}{z}} \}$ and $\pmast 20\pmcdot1$ we
  have:
  \begin{flalign*}
    \quad \pmSub{\pmast 20\pmcdot1}
          {\{\pmcls{z}{ \pmpred{\chi}{z} }  = \pmcls{z}{\pmpf{\phi}{z}} \}}
          {f\{\pmcls{z}{\pmpf{\varphi}{z}}\}}\
          \pmthm \pmdottt
          \{\pmcls{z}{ \pmpred{\chi}{z} }  = \pmcls{z}{\pmpf{\phi}{z}} \}
          \pmdot \pmiff \pmdott & \\
          \pmsome{\psi} \pmdott \pmpred{\psi}{z} \pmdot
          \pmiff_{z} \pmdot \pmpf{\phi}{z}
          \pmdott
          \{\pmcls{z}{ \pmpred{\chi}{z} }  = \pmpred{\psi}{\pmhat{z}} \} \label{equation_chwistek_5} && \tag{6} 
  \end{flalign*}
  Hence, if we take a restricted view of the scope of $\pmcls{z}{\pmpf{\phi}{z}}$ in Equation \ref{equation_chwistek_4}, we might try to apply
  $\pmast 20\pmcdot01$ to thus eliminate a defined expression from $\pmnot \{\pmcls{z}{ \pmpred{\chi}{z} }  = \pmcls{z}{ \pmpf{\phi}{z} }   \}$
  in Equation \ref{equation_chwistek_4}, so as to obtain Equation \ref{equation_chwistek_3}:\footnote{The outermost braces on the left of
  Equation \ref{equation_chwistek_3} are of course unnecessary in light of Definition $\pmast 9\pmcdot021$. Since intuitions apparently vary regarding the truth or
  falsity of such equations it may be that some redundant punctuation is helpful.}
  \begin{flalign*}
    \quad
    \pmthm \pmdottt
    \pmnot \{ \pmsome{\psi} \pmdott \pmpred{\psi}{z} \pmdot
    \pmiff_{z} \pmdot  \pmpf{\phi}{z}  \pmdott
    \{\pmcls{z}{ \pmpred{\chi}{z} }  = \pmpred{\psi}{\pmhat{z}} \} \}
    \pmdott \pmiff \pmdott & \\
    \pmsome{\psi} \pmdott \pmpred{\psi}{z} \pmdot \pmiff_{z} \pmdot
    \pmpf{\phi}{z}  \pmandd
    \pmnot \{\pmcls{z}{ \pmpred{\chi}{z} }  = \pmpred{\psi}{\pmhat{z}} \} \label{equation_chwistek_3} && \tag{7}
  \end{flalign*}
  If however we take the scope for $\pmcls{z}{ \pmpf{\phi}{z} }$ as the whole of $\pmnot \{\pmcls{z}{ \pmpred{\chi}{z} }  = \pmcls{z}{ \pmpf{\phi}{z} }   \}$
  we obtain the following in using $\pmast 20\pmcdot01$ to eliminate a defined expression from Equation \ref{equation_chwistek_4}:
  \begin{flalign*}
    \quad
    \pmthm \pmdottt
    \pmsome{\psi} \pmdott \pmpred{\psi}{z}
    \pmdot \pmiff_{z} \pmdot
    \pmpf{\phi}{z} \pmdott
    \pmnot \{ \pmcls{z}{ \pmpred{\chi}{z} }  = \pmpred{\psi}{\pmhat{z}} \} 
    \pmdott \pmiff \pmdott & \\
    \pmsome{\psi} \pmdott \pmpred{\psi}{z}
    \pmdot \pmiff_{z} \pmdot
    \pmpf{\phi}{z} \pmandd
    \pmnot \{\pmcls{z}{ \pmpred{\chi}{z} }  = \pmpred{\psi}{\pmhat{z}} \} \label{equation_chwistek_3_1} && \tag{7.1}
  \end{flalign*}
  As Chwistek aims to show that Equation \ref{equation_chwistek_3} yields pathological results let's consider this Equation before considering
  whether the scope of $\pmcls{z}{ \pmpf{\phi}{z} }$ can be taken as less than the whole of $\pmnot \{\pmcls{z}{ \pmpred{\chi}{z} }  = \pmcls{z}{ \pmpf{\phi}{z} }   \}$.
  At this point the first hurdle for Chwistek's argument arises. 
  \begin{lemma}\label{lemma_equation_7}
    On the basis of an informal interpretation, Equation \ref{equation_chwistek_3} does not imply any flaw in \emph{Principia}'s definition of a class,
    provided that, in the context of use, extensionally equivalent predicative functions are equal.
  \end{lemma}
  \begin{proof}
    The validity of Equation \ref{equation_chwistek_3} may be demonstrated informally using an interpretation of the asserted propositional function.
    In informal English, Equation \ref{equation_chwistek_3} states that:
    \begin{interpretation}[Informal statement of Equation \ref{equation_chwistek_3}]\label{interpertation_equation_seven}
      For a definite propositional function $\pmpf{\phi}{z}$ and $\pmpred{\chi}{z}$ any predicative function such
      that $\pmpred{\chi}{z} \pmdot \pmiff_{z} \pmdot \pmpf{\phi}{z}$ is significant /
      a proposition: the first of the following biconditionals is true if and only if the second of the following biconditionals is true.
      \begin{enumerate}
      \item (Left-Hand Side of Equation \ref{equation_chwistek_3}) There does not exist a predicative function
        $\pmpred{\psi}{z}$ which is (simultaneously) formally equivalent to both of the functions $\pmpf{\phi}{z}$ and $\pmpred{\chi}{z}$.
      \item (Right-Hand Side of Equation \ref{equation_chwistek_3}) There exists a predicative function
        $\pmpred{\psi}{z}$ which is both: (i) formally equivalent to $\pmpf{\phi}{z}$ and (ii) \emph{not} formally equivalent to the predicative function $\pmpred{\chi}{z}$.
      \end{enumerate}
    \end{interpretation}
    To provide a precise interpretation of Equation \ref{equation_chwistek_3} we need to define the notions of ``Truth'' and ``Falsity'', of each type, for \emph{Principia}.
    This is known to be a difficult problem (\cite{tarski1936}: \S 5). Nevertheless, if we accept \emph{Principia}'s primitive, typed notions of Truth and Falsity as
    given, then with the eye of faith, one can see using Interpretation \ref{interpertation_equation_seven} that Equation \ref{equation_chwistek_3} is valid by considering two
    classes of cases that form a partition of the set of all possible cases to consider: $\pmpred{\chi}{z} \pmdot \pmiff_{z} \pmdot \pmpf{\phi}{z}$ and
    $\pmnot \{\pmpred{\chi}{z} \pmdot \pmiff_{z} \pmdot \pmpf{\phi}{z}\}$.

    Before entering into details, it will simplify the discussion to consider
    the Right-Hand Side of Equation \ref{equation_chwistek_3} with some defined expressions eliminated via a second application of $\pmast 20\pmcdot01$
    (Cf. \cite{pm1910v1} proof of $\pmast 20\pmcdot13$):
    \begin{flalign*}
      \quad \pmthm \pmdottt
      \pmnot \{ \pmsome{\psi} \pmdott & \pmpred{\psi}{z} \pmdot
      \pmiff_{z} \pmdot  \pmpf{\phi}{z}  \pmdott
      \{\pmcls{z}{ \pmpred{\chi}{z} }  = \pmpred{\psi}{\pmhat{z}} \} \}
      \pmdott \pmiff \pmdott & \\
      & \pmsome{\psi, \Psi} \pmdott \pmpred{\psi}{z} \pmdot \pmiff_{z} \pmdot
      \pmpf{\phi}{z}  \pmandd
      \pmpred{\Psi}{x} \pmdot \pmiff_x \pmdot \pmpred{\chi}{x} \pmandd  \pmnot \{\pmpred{\Psi}{\pmhat{z}}  = \pmpred{\psi}{\pmhat{z}} \}  
      \tag{  \ref{lemma_equation_7}.0} \label{equation_chwistek_5.0} &
    \end{flalign*}
    \begin{description}
    \item[Case $\pmpred{\chi}{z} \pmdot \pmiff_{z} \pmdot \pmpf{\phi}{z}$] 
      The case hypothesis implies that both sides of Equation \ref{equation_chwistek_3} are false:
      \begin{enumerate}
      \item The left-hand side of Equation \ref{equation_chwistek_3} is false since, as $\pmpred{\chi}{z}$ \emph{is}, by hypothesis, formally equivalent to $\pmpf{\phi}{z}$,
        there does exist a predicative function
        $\pmpred{\psi}{z}$ which is simultaneously formally equivalent to both $\pmpred{\chi}{z}$ and $\pmpf{\phi}{z}$ ($\pmpred{\chi}{z}$ itself being just such a function).
      \item To show that the right-hand side of Equation \ref{equation_chwistek_3} is false, assume to derive a contradiction that it is true -
        that there exists the required predicative function $\pmpred{\psi}{z}$. By examining Equation \ref{equation_chwistek_5.0}, we see that
        (i) $\pmpred{\psi}{z}$ must be formally equivalent to $\pmpf{\phi}{x}$, and (ii) there must also exist a predicative function $\pmpred{\Psi}{x}$
        that is formally equivalent to $\pmpred{\chi}{x}$ and yet (iii)  $\pmpred{\psi}{z}$ and $\pmpred{\Psi}{x}$ are not equal.
        But by the case hypothesis that is impossible, since the formal equivalence of $\pmpred{\chi}{z}$ and $\pmpf{\phi}{z}$ implies - given (i), (ii),
        and the assumption that, in the context of use, extensionally equivalent predicative functions are equal - 
        that $\pmpred{\psi}{z}$ and $\pmpred{\Psi}{x}$ \emph{are} formally equivalent and not equal, a contradiction. Hence there does \emph{not} exists
        any such predicative function and the right-hand side of Equation \ref{equation_chwistek_3} is false as required.
      \end{enumerate}
    \item[Case $\pmnot \{\pmpred{\chi}{z} \pmdot \pmiff_{z} \pmdot \pmpf{\phi}{z}\}$] 
      The case hypothesis implies that both sides of Equation \ref{equation_chwistek_3} are true:
      \begin{enumerate}
      \item The left-hand side of Equation \ref{equation_chwistek_3} is true since as $\pmpred{\chi}{z}$ is \emph{not}, by hypothesis, formally equivalent to $\pmpf{\phi}{z}$,
        there does not exist any predicative function
        $\pmpred{\psi}{z}$ which is simultaneously formally equivalent to both $\pmpred{\chi}{z}$ and $\pmpf{\phi}{z}$.
      \item The right-hand side of Equation \ref{equation_chwistek_3} is true since: firstly, by an appropriate instance of the axiom of reducibility, there does exist a predicative function
        $\pmpred{\psi}{z}$ which is formally equivalent to $\pmpf{\phi}{z}$; secondly, as $\pmpred{\chi}{z}$ is \emph{not}, by hypothesis, formally equivalent to $\pmpf{\phi}{z}$,
        the aforementioned $\pmpred{\psi}{z}$ is \emph{not} formally equivalent to $\pmpred{\chi}{z}$.
      \end{enumerate}
    \end{description}
  \end{proof}
  A comment on Russell and Whitehead's (fixed) interpretation of the formal system of \emph{Principia} may assist in clarifying the context restriction used in Lemma \ref{lemma_equation_7}.
  The key point to note is that \emph{Principia}'s variables in an asserted propositional function concern any - or,  in certain propositional functions, some - entity of the relevant type:
  \begin{remark}\label{remark_context}
    The hypothesis that ``in the context of use, extensionally equivalent predicative functions are equal'' excludes very few contexts encountered in \emph{Principia}:
    \begin{enumerate}
    \item For predicative functions of an individual only, functions formed using (i) names of specific constants (for individuals and or relations / properties or individuals)
      and / or, (ii) Use of symbols other than the standard logical constants ($\pmnot$, $\pmor$, $\pmpred{\phi}{x}$ etc.)
      and associated symbols in forming the (interpreted) symbol-sequence for the function (in addition to symbols for variables for objects of a given type) are
      the only excluded functions.
    \item Likewise, for predicative functions of a higher type the only excluded functions involve either (i) specific constants naming members of a given type; and / or,
      (ii) Use of symbols other than the standard logical constants ($\pmnot$, $\pmor$, $\pmpred{\phi}{x}$ etc.)
      and associated symbols in forming the (interpreted) symbol-sequence for the function (in addition to symbols for variables for objects of a given type).
    \end{enumerate}
  \end{remark}
  A formal proof or refutation of Equation \ref{equation_chwistek_3} must involve a representation of Equation \ref{equation_chwistek_3} in some formal system and hence
  assume the validity of that system.
  As Chwistek's aim is, effectively, to challenge the consistency of the formal system of \emph{Principia} a demonstration that uses this system
  begs the question of whether this system is consistent.\footnote{If a false proposition / propositional function of some type is provable in \emph{Principia} then
  any proposition  / propositional function of that type is provable in the system using modus ponens.} Nevertheless, since Chwistek's claim that
  Equation \ref{equation_chwistek_3} is false has attracted some support it may be useful to
  consider a second perspective on whether Equation \ref{equation_chwistek_3} is a theorem of \emph{Principia}.

  \begin{lemma}\label{lemma_equation_7_1}
    By \emph{Principia}'s lights, Equation \ref{equation_chwistek_3} is a theorem.
  \end{lemma}
  \begin{proof}
    By an instance of $\pmast 20\pmcdot15$:\footnote{To anticipate an issue concerning \emph{Principia}'s substitution rules discussed below, the reader should note
    that $\pmpred{\chi}{z}$ is not an existing free or real variable of $\pmast 20\pmcdot15$ prior to this substitution. The resemblance between
    `` $\pmpred{\chi}{z}$'' and ``$\pmpf{\chi}{z}$'' in this context is a mere typographic accident, as in substance these are not the same function.}
    \begin{flalign*}
      \quad \pmSubb{\pmast 20\pmcdot15}
            {\pmpred{\chi}{z} }
            {\pmpf{\psi}{z}}
            {\pmpf{\phi}{z}}
            {\pmpf{\chi}{z}}\
            \pmthm \pmdottt
            \pmpred{\chi}{z} \pmdot \pmiff_{z} \pmdot \pmpf{\phi}{z}
            \pmdott \pmiff \pmdot              
            \pmcls{z}{ \pmpred{\chi}{z} }  = \pmcls{z}{ \pmpf{\phi}{z} }   
            \tag{  \ref{lemma_equation_7}.1} \label{equation_chwistek_5.1} &&
    \end{flalign*}
    Hence, if we use $\pmast 20\pmcdot01$ to eliminate a defined expression from $\pmcls{z}{ \pmpred{\chi}{z} }  = \pmcls{z}{ \pmpf{\phi}{z} }$
    in Equation \ref{equation_chwistek_5.1}, we obtain:
    \begin{flalign*}
      \quad [(\ref{equation_chwistek_5.1}) \pmand \pmast 20\pmcdot01]\
      \pmthm \pmdottt\
      & \pmpred{\chi}{z} \pmdot \pmiff_{z} \pmdot \pmpf{\phi}{z} 
      \pmdott \pmiff \pmdott & \\
      & \pmsome{\psi} \pmdott \pmpred{\psi}{z} \pmdot
      \pmiff_{z} \pmdot
      \pmpf{\phi}{z}
      \pmdott
      \{\pmcls{z}{ \pmpred{\chi}{z} }  = \pmpred{\psi}{\pmhat{z}} \} 
      \tag{  \ref{lemma_equation_7}.2} \label{equation_chwistek_5.2} &
    \end{flalign*}
    We also have however:
    \begin{flalign*}
      \quad [(\ref{equation_chwistek_5.1}) \pmand \pmast 4\pmcdot11]\
      \pmthm \pmdottt\
      \pmnot \{\pmpred{\chi}{z} \pmdot \pmiff_{z} \pmdot \pmpf{\phi}{z} \}
      \pmdott \pmiff \pmdot              
      \pmnot \{ \pmcls{z}{ \pmpred{\chi}{z} }  = \pmcls{z}{ \pmpf{\phi}{z} } \}
      \tag{  \ref{lemma_equation_7}.3} \label{equation_chwistek_5.3} &
    \end{flalign*}
    As before, using $\pmast 20\pmcdot01$ to eliminate a defined expression from $\pmnot \{\pmcls{z}{ \pmpred{\chi}{z} }  = \pmcls{z}{ \pmpf{\phi}{z} }   \}$
    in Equation \ref{equation_chwistek_5.3}, we obtain:
    \begin{flalign*}
      \quad [(\ref{equation_chwistek_5.3}) \pmand \pmast 20\pmcdot01]\
      \pmthm \pmdottt\
      & \pmnot \{\pmpred{\chi}{z} \pmdot \pmiff_{z} \pmdot \pmpf{\phi}{z} \} 
      \pmdott \pmiff \pmdott  & \\
      & \pmsome{\psi} \pmdott \pmpred{\psi}{z}
      \pmdot \pmiff_{z} \pmdot
      \pmpf{\phi}{z}
      \pmdott
      \pmnot \{\pmcls{z}{ \pmpred{\chi}{z} }  = \pmpred{\psi}{\pmhat{z}} \}
      \tag{  \ref{lemma_equation_7}.4} \label{equation_chwistek_5.4} &
    \end{flalign*}
    An abbreviated summary of the remaining steps for a proof that Equation \ref{equation_chwistek_3} is a theorem of \emph{Principia} is as follows.
    \begin{description}
    \item[Case $\pmpred{\chi}{z} \pmdot \pmiff_{z} \pmdot \pmpf{\phi}{z}$] 
      The case hypothesis, Equation \ref{equation_chwistek_5.2},  and $\pmast 4\pmcdot11$ imply that
      the left-hand side of Equation \ref{equation_chwistek_3} is false. Similarly, the case hypothesis, Equation \ref{equation_chwistek_5.4}, and $\pmast 4\pmcdot11$ imply that
      the right-hand side of Equation \ref{equation_chwistek_3} is false. 
    \item[Case $\pmnot \{\pmpred{\chi}{z} \pmdot \pmiff_{z} \pmdot \pmpf{\phi}{z}\}$] 
      The case hypothesis, Equation \ref{equation_chwistek_5.2},  and $\pmast 4\pmcdot11$ imply that
      the left-hand side of Equation \ref{equation_chwistek_3} is true. Similarly, the case hypothesis, Equation \ref{equation_chwistek_5.4}, and $\pmast 4\pmcdot11$ imply that
      the right-hand side of Equation \ref{equation_chwistek_3} is true. 
    \end{description}
  \end{proof}

  With an additional application of $\pmast 20\pmcdot01$
  to Equation \ref{equation_chwistek_5.2}, we can use the function $\pmcls{z}{ \pmpred{\chi}{z} }  = \pmcls{z}{ \pmpf{\phi}{z} }$ to illustrate 
  \emph{Principia}'s claim that ``$f\{\pmcls{z}{\pmpf{\psi}{z}}\}$ is always an \emph{extensional} function of $\pmpf{\psi}{\pmhat{z}}$'' (\cite{pm1910v1}: 197).
  \begin{flalign*}
    \quad [(\ref{equation_chwistek_5.2}) \pmand \pmast 20\pmcdot1]\
    \pmthm \pmdottt\
    & \pmpred{\chi}{z} \pmdot \pmiff_{z} \pmdot \pmpf{\phi}{z} 
    \pmdott \pmiff \pmdott & \\
    & \pmsome{\psi, \Psi} \pmdott \pmpred{\psi}{z} \pmdot \pmiff_{z} \pmdot
    \pmpf{\phi}{z}  \pmandd
    \pmpred{\Psi}{x} \pmdot \pmiff_x \pmdot \pmpred{\chi}{x} \pmandd \{\pmpred{\Psi}{\pmhat{z}}  = \pmpred{\psi}{\pmhat{z}} \} 
    \tag{  \ref{lemma_equation_7}.5} \label{equation_chwistek_5.5} &
  \end{flalign*}
  For reasons set out at Remark \ref{remark_context}, the context of Equation \ref{equation_chwistek_5.5} is such that predicative propositional functions
  of a given type and arity are equal if and only if they are true for the same arguments.
  
  As noted above, a demonstration that Equation \ref{equation_chwistek_3} is valid which uses \emph{Principia}'s definition of a class simply begs the question of whether
  \emph{Principia}'s definition leads to paradox, since Equation \ref{equation_chwistek_3} is derived from this very same approach. We should
  therefore consider Chwistek's reasons for rejecting Equation \ref{equation_chwistek_2}. Chwistek's terse conclusion that Equation \ref{equation_chwistek_2} ``is evidently false''
  (\cite{chwistek1922}: 339) is elaborated upon elsewhere
  (\cite{chwistek1924} in \cite{Linsky2004}). In this source,
  Chwistek considers a slightly different function $f$ to that mentioned above:
  \begin{flalign*}
    f\{\pmcls{z}{\pmpf{\psi}{z}}\} = \pmnot \{\pmcls{z}{\pmpf{\psi}{z}}  = \pmpred{\theta}{\pmhat{z}} \}\ \label{item_chwistek_item_3} \text{and} & \tag{8} \\
    f\{\pmpred{\Phi}{\pmhat{z}}\} = \pmnot \{ \pmpred{\Phi}{\pmhat{z}} = \pmpred{\theta}{\pmhat{z}} \}\ \label{item_chwistek_item_4} \text{etc.} && \tag{9}
  \end{flalign*}
  
  Using this function $f$ and different names
  for variables Chwistek considers, instead of Equation \ref{equation_chwistek_2}, the following Equation:\footnote{\cite{chwistek1924} in \cite{Linsky2004}: 68,
  modified through changing of the reference number for this equation SB. The preceding equation in the derivation
  has ``$\pmsome{\psi}$'', but this should actually be ``$\pmsome{\Phi}$''.}
  \begin{flalign*}
    \quad \pmnot \pmsome{\Phi} \pmdott \pmpred{\Phi}{x}
    \pmdot \pmiff_{x} \pmdot
    \pmpf{\psi}{x} \pmdott
    \{\pmpred{\Phi}{\pmhat{z}}  = \pmpred{\theta}{\pmhat{z}} \}
    \pmdott \pmiff \pmdott
    \pmsome{\Phi} \pmdott \pmpred{\Phi}{x} \pmdot
    \pmiff_{x} \pmdot
    \pmpf{\psi}{x}\ \pmandd & \nonumber \\
    \pmnot \{\pmpred{\Phi}{\pmhat{z}}   = \pmpred{\theta}{\pmhat{z}} \}. \label{equation_chwistek_6} \tag{10} &&
  \end{flalign*}
  In relation to Equation \ref{equation_chwistek_6} Chwistek then proposes ``Now let us take $\pmpred{\theta}{\pmhat{z}}$ for $\pmpf{\psi}{\pmhat{z}}$''
  (\cite{chwistek1924} in \cite{Linsky2004}: 68). The
  result of such a substitution is the following:
  \begin{flalign*}
    \quad \pmnot \pmsome{\Phi} \pmdott \pmpred{\Phi}{x}
    \pmdot \pmiff_{x} \pmdot
    \pmpred{\theta}{x} \pmdott
    \{\pmpred{\Phi}{\pmhat{z}}  = \pmpred{\theta}{\pmhat{z}} \}
    \pmdott \pmiff \pmdott
    \pmsome{\Phi} \pmdott \pmpred{\Phi}{x} \pmdot
    \pmiff_{x} \pmdot
    \pmpred{\theta}{x}\ \pmandd & \nonumber \\
    \pmnot \{\pmpred{\Phi}{\pmhat{z}}   = \pmpred{\theta}{\pmhat{z}} \}. \label{equation_chwistek_7} && \tag{11}
  \end{flalign*}
  
  Surprisingly, neither Chwistek nor Linsky discuss whether the substitution required to derive Equation \ref{equation_chwistek_7} is reasonable.
  Whilst \emph{Principia} does not present substitution
  rules that explicitly prohibit the substitution of $\pmpred{\theta}{\pmhat{z}}$ for $\pmpf{\psi}{\pmhat{z}}$ in Equation \ref{equation_chwistek_6}, this
  is an \emph{error of omission}.
  \begin{lemma}\label{lemma_equation_7_2}
    The substitution used to derive Equation \ref{equation_chwistek_7} from Equation \ref{equation_chwistek_6}
    is not valid in the sense that it permits the inference of false propositional functions from true.
  \end{lemma}
  \begin{proof}
    A comprehensive statement of \emph{Principia}'s substitution rules is beyond the scope of this paper. In lieu of such, the problem with the
    class of substitutions in question can be illustrated by a comparison with first-order arithmetic. That is, it can be seen that
    well-formulated substitution rules for \emph{Principia} would prohibit such a substitution, for reasons illustrated in the following spurious demonstration
    that $(x_1 < x_1)$ is a theorem of first-order arithmetic (PA):\footnote{The ``first-order arithmetic'' discussed here is of course a metamathematical theory which
    has no real connection with \emph{Principia}'s arithmetic. An example from an entirely different kind of formal system is nevertheless appropriate since the point is to
    illustrate an issue relating to substitution rules that affects both kinds of systems.} \\

    \noindent \begin{minipage}{\textwidth}
      \centering
      \textit{Spurious demonstration of $\vdash_{PA} (x_1 < x_1)$}
      \medskip 
      \begin{tabular}[htbp]{ l l l }
        (1) & $\vdash_{PA} (\forall x_1)(\exists x_2) (x_1 < x_2)$ & Established theorem\\
        (2) & $\vdash_{PA} (\exists x_2) (x_1 < x_2)$ & (1), universal instantiation, $x_1$ free for $x_1$\\
        (3) & $\vdash_{PA}  (x_1 < c)$ & (2), Rule C / existential instantiation,\\
        &           & $c$ a new constant\\
        (4) & $\vdash_{PA}  (x_1 < x_1)$          & (3), by an \emph{invalid} substitution rule\
      \end{tabular}
    \end{minipage}
  \end{proof}
  The spurious demonstration follows Mendelson's example of a fallacious inference allowed by a misuse of Rule C wherein the restriction concerning the use of Generalisation
  is ignored (\cite{mendelson2015}: \S 2.6). Whilst much of the spurious demonstration is a paraphrase of Mendelson's example an important difference to note is that the former
  does \emph{not} use Generalisation. Mendelson's version of Rule C is thus silent with respect to such spurious proofs - the substitution of ``$x_1$'' for ``$c$'' in passing from
  (3) to (4) is officially neither prohibited nor approved.
  
  In comparing the spurious PA demonstration with Chwistek's substitution of $\pmpred{\theta}{\pmhat{z}}$ for $\pmpf{\psi}{\pmhat{z}}$ in Equation \ref{equation_chwistek_6},
  I treat $\pmpred{\theta}{\pmhat{z}}$ as a variable and $\pmpf{\psi}{\pmhat{z}}$ as a constant expression.
  \emph{Principia}'s variables, after $\pmast 12$, only range over individuals and matrices (\cite{pm1910v1}: 172). Since $\pmpred{\theta}{\pmhat{z}}$ is
  a predicative function, it may be viewed as variable, comparable to $x_1$, whereas $\pmpf{\psi}{\pmhat{z}}$ is comparable to a constant expression, like $c$. In the terminology
  of the First Edition of \emph{Principia}, in the context of an asserted propositional function (such as an instance of $\pmast 20\pmcdot1$):
  \begin{enumerate}
  \item Where $\pmpred{\theta}{\pmhat{z}}$ is a real variable, it is ``any'' predicative function of the appropriate type, i.e. it is
    an arbitrary predicative function or, in other words, a specific predicative function of the appropriate type but it is not determined which;
  \item $\pmpf{\psi}{\pmhat{z}}$ is a specific propositional function of a definite type.
  \end{enumerate}

  At this point we may return to the question of the scope of class symbols, such as ``$\pmcls{z}{ \pmpf{\phi}{z} }$'' in applying  $\pmast 20\pmcdot01$.
  \begin{lemma}\label{lemma_equation_7_3}
    Chwistek's derivation of Equation \ref{equation_chwistek_2} from Equation \ref{equation_chwistek_1} exhibits an error in suggesting that,
    in general, within a given context,
    two distinct determinations of the scope of
    $\pmcls{z}{\pmpf{\phi}{z}}$ may be used in applying $\pmast 20\pmcdot01$ to eliminate a defined term.
  \end{lemma}
  \begin{proof}
    \emph{Principia}'s discussion of the scope of descriptions illustrates how different determinations of the scope of the definiendum may, in some contexts, result
    in different, potentially inconsistent, propositions being determined to be the associated definiens (\cite{pm1910v1}: 181-2). Since Equation \ref{equation_chwistek_3}
    is a theorem, the two distinct choices that Chwistek makes in the derivation of Equation \ref{equation_chwistek_2} do not, in this context, result in an inconsistent
    proposition being established as a theorem. Nevertheless, Chwistek's equivocation concerning the scope of the defined term to be eliminated is a deplorable
    misrepresentation of the system of \emph{Principia}. Where a second application of $\pmast 20\pmcdot01$ to an equation such as Equation \ref{equation_chwistek_4} is
    required, a consistent determination of the scope of a definiendum is clearly required. On this basis, we can easily determine that Equation \ref{equation_chwistek_3_1},
    rather than Equation \ref{equation_chwistek_3}, is the correct result of applying $\pmast 20\pmcdot01$ to Equation \ref{equation_chwistek_4}.
    As indicated at Remark \ref{remark_context_0}, choosing different scopes in most contexts encountered in \emph{Principia} will make no difference. As one of Chwistek's
    aims however is to challenge \emph{Principia}'s rule for determining when propositional functions are equal (\cite{chwistek1924} in \cite{Linsky2004}: 69),
    Chwistek's equivocation on this point is untenable.
  \end{proof}
  It is somewhat ironic that Chwistek's misrepresentations of the system of \emph{Principia} in some places result in theorems of \emph{Principia}
  being derived via demonstrations that use rules or conventions that are not part of the system of \emph{Principia}. If the propositions / propositional functions
  involved were not, coincidentally, theorems of \emph{Principia} one could simply point out the disparity with respect to the convention or rule of substitution used.
  Where theorems are involved, it has been important to consider Chwistek's reasons for claiming that these are false.
\end{proof}
\section{Discussion}
Chwistek's fumbles over scope rules and rules of substitution
does highlight a lack of clear guidance on this in \emph{Principia}. Russell and Whitehead are hardly however the only logicians from this period who neglected to provide a comprehensive, correct
set of such statements for their system (\cite{pager1962}). The rather uncritical subsequent reception of Chwistek's claim possibly reflects a bias towards the view that there is ``something rotten''
in the foundations of \emph{Principia}, though the dissemination of acceptance of Chwistek's claim may partly be to blame for the emergence of this view. Russell himself seemed hesitant to accept
the detail of Chwistek's demonstration, though as he took it that a contradiction was being deduced from the axiom of reducibility, which he was no longer committed
to defending, he did not regard the claim as a cause for concern (\cite{Linsky2004}: 64).

G\"odel's erroneous comments on this point have no doubt contributed to an unfair sullying of \emph{Principia}'s reputation:
\begin{quote}
If this latter rule [``of replacing defined symbols by their definiens''] is applied to expressions containing other defined symbols it re-
quires that the order of elimination of these be indifferent. This however is
by no means always the case ($\pmpred{\phi}{\pmhat{u}} = \pmcls{u}{\pmpf{\phi}{u}}$, e.g., is a counter-example). (\cite{godel1944}: 120,
modified through interpolation of quote from previous sentence in square brackets)
\end{quote}
It is very easy to demonstrate that the rule of thus replacing defined symbols by their definiens does not require that the order of elimination be indifferent.
\emph{Principia}'s definitions are introduced in a definite order which, putting aside for the moment the type of the symbols involved, is explicitly provided by Russell and Whitehead's
system for numbering Primitive propositions, definitions, and theorems. The symbols for descriptions, for example,  $\pmast 14\pmcdot01$ etc. are defined subsequent to the
symbols for equality ($\pmast 13\pmcdot01$ etc.). From a logical point of view, the elimination of defined symbols from a proposition or propositional function of \emph{Principia},
other things being equal, ought to proceed in reverse order to the order in which definitions are introduced. Thus, logically speaking, symbols for descriptions ought to be eliminated prior
to eliminating symbols for equality. G\"odel's nonsense remarks set up an opponent of straw that simply distracts attention from the genuine task of trying to
understand how a system such as \emph{Principia} actually works.

\section{Conclusion}

Contemporary logicians overwhelmingly view \emph{Principia} as of historical interest at best. Approaches to defining the syntax and semantics of formal systems pioneered by G\"odel, Tarski
and others from the 1920's onwards provide, on this account, a better way of doing things.

Having slipped from the limelight, \emph{Principia} is at risk of being misrepresented
due to erroneous, historical criticisms failing to solicit the interest required to generate a rebuttal. The challenge of examining received accounts of \emph{Principia}
is compounded by the fact that contemporary scholars are, for the most part, reluctant to accept Whitehead and Russell's definition of
the formal system of \emph{Principia}. To focus on the issue raised by Chwistek's argument the discussion above skirts around a variety of complex
questions concerning the contents of \emph{Principia}, such as the status of \emph{Principia}'s propositions / propositional functions, it's Primitive propositions (axioms and
inference rules) and so on. Whilst clearly these are important issues, discussion of such matters should not get in the way of calling out false criticism. Whatever one's broader views on
these things, it is clear Chwistek's claim that with \emph{Principia}'s definition of a  class Richard’s paradox can be formulated is false. The faulty inferences connected with Chwistek's
demonstration are due to incorrect uses of scope and substitution rule that is no part of \emph{Principia}; they are not, as Chwistek hypothesised, evidence of an error in \emph{Principia}'s
definition of a class.

\end{document}